\documentclass[12pt]{amsart}

\usepackage{amsmath}	
\usepackage{amssymb}

\newtheorem{thm}{Theorem}[section]
\newtheorem*{thmA}{Theorem A}
\newtheorem*{thmB}{Theorem B}
\newtheorem{lemma}[thm]{Lemma}
\newtheorem{prop}[thm]{Proposition}
\newtheorem{cor}[thm]{Corollary}

\newtheorem{example}[thm]{Example}

\theoremstyle{remark}

\numberwithin{equation}{section}

\newcommand{\bN}{\mathbb{N}}

\newcommand{\bQ}{\mathbb{Q}}
\newcommand{\bR}{\mathbb{R}}
\newcommand{\bZ}{\mathbb{Z}}

\newcommand{\cC}{{\mathcal{C}}}

\newcommand{\ub}{{\mathbf{b}}}
\newcommand{\uu}{{\mathbf{u}}}
\newcommand{\uv}{{\mathbf{v}}}
\newcommand{\ux}{{\mathbf{x}}}
\newcommand{\uy}{{\mathbf{y}}}
\newcommand{\uzero}{{\mathbf{0}}}

\newcommand{\Aff}{\mathrm{Aff}\,}
\newcommand{\et}{\quad\text{and}\quad}

\newcommand{\Inf}{\mathrm{Inf}\,}

\newcommand{\proj}{\mathrm{proj}}
\newcommand{\Set}[2]{\left\{#1 \vphantom{#2}\,\left|
            \vphantom{#1}\,\,#2\right. \right\}}
\newcommand{\un}{\mathbf{1}}

\begin{document}

\baselineskip=15.2pt 

\date{\today}

\title[One-sided approximation]
{One-sided approximation in affine function spaces}
\author{David Handelman and Damien Roy}

\thanks{Both authors supported in part by NSERC Discovery grants.}

\begin{abstract}
Let $H$ be a subgroup of a partially ordered abelian
group $G$ with order unit $u$, and let $S(G,u)$ denote
the convex subset of $\bR^G$ consisting of all traces
(states) $\tau$ on $G$ with $\tau(u)=1$.  We say that
$H$ has property $(B)$ if, for any integer $m\ge 2$,
any $h\in H$ and any $\epsilon>0$, there exists
$h'\in H$ such that $\tau(h)-m\tau(h')\ge -\epsilon$ for each
$\tau\in S(G,u)$.  We show that, if $S(G,u)$ is
finite-dimensional, this condition is equivalent to asking that
$\tau(H)$ is $\{0\}$ or dense in $\bR$ for all $\tau$
in the smallest face of $S(G,u)$ containing all traces
that vanish identically on $H$.  When $G$ is a simple
dimension group and $H$ is a convex subgroup of $G$, we
show that $G/H$ is unperforated if and only if $H$ has
property $(B)$.  We apply both results to provide a
criterion for a trace of $G$ to be refinable when $G$
is a simple dimension group with finitely many pure traces.
\end{abstract}

\maketitle

\section{Introduction}
\label{sec:intro}

Throughout, $n$ will mean a positive integer.
We denote by $(\bR^n)^*$, the dual space to $\bR^n$, consisting
of the linear functionals on $\bR^n$; by $(\bR^n)^{+}$, the set of
elements of $\bR^n$ with non-negative coordinates; and
by $(\bR^n)^{++}$, the set of those with strictly positive
coordinates.  We also denote by $\|x\|$ the Euclidean
norm of a point $x\in\bR^n$.

Let $H$ be a subgroup of $\bR^n$, and let $m>1$ be an integer.
The condition
\begin{itemize}
\item[$(A_m)$\quad]
      for all $\epsilon > 0$ and $h \in H$, there exists
      $h'\in H$ such that $\|h - mh' \| \le \epsilon$
\end{itemize}
is independent of $m$ and is equivalent to $H$ being
dense in the vector space $\bR\,H$ that it spans; so we
may as well refer to it as property $(A)$. By a theorem
of Kronecker, this in turn implies the following.

\begin{thmA}
Let $m$ and $H$ be as above.  Then $H$ satisfies
property $(A_m)$ if and only if for all $\tau \in (\bR^n)^*$,
either $\tau(H) = \{0\}$ or $\tau(H)$ is dense in $\bR$.
\end{thmA}

Here, we are interested in the following one-sided approximation
property:
\begin{itemize}
\item[$(B_m)$\quad]
 for all $\epsilon > 0$ and $h \in H$, there exists $h' \in H$
 such that all coordinates of $h - mh'$ are bounded below by $-\epsilon$.
\end{itemize}
\noindent
For example, the discrete subgroup $H = \bZ^n$ of $\bR^n$ has this property,
while $H= \bZ(1,-1)$ inside $\bR^2$ does not.

Our principal result, in this context, characterizes subgroups $H$
satisfying this property, in terms of \emph{positive} linear functionals
in a fashion analogous to that in Theorem A.  Recall that a positive linear
functional on $\bR^n$ is a linear functional sending $(\bR^n)^+$ to
$\bR^+$ or, equivalently, that sends the standard basis elements to nonnegative real
numbers. Let $K_n$ denote the usual standard $(n-1)$-simplex consisting
of positive linear functionals $\tau$ in $(\bR^n)^*$ such that
$\tau(1,\dots, 1) = 1$;
its vertices are the $n$ coordinate functions $\tau_1,\dots,\tau_n$ where
$\tau_i\colon\bR^n\to \bR$ is projection onto the $i$th coordinate.

\begin{thmB}
Let $m$ and $H$ be as above. Define $F$ to be the
smallest face of $K_n$ containing the set
\[
 Z(H):=\Set{\tau \in K_n}{\tau(H) = \{0\}}.
\]
Then $H$ satisfies property $(B_m)$ if and only if for all
$\tau \in F$, either $\tau(H) = \{0\}$ or $\tau(H)$ is dense in $\bR$.
\end{thmB}

\medskip
The set $Z(H)$ is a compact convex subset of $K_n$.
If it is empty, then $F$ is empty and the condition is vacuous,
so $H$ satisfies $(B_m)$ in this case.  We will see in section
\ref{sec:suf} that this happens if and only if $H\cap(\bR^n)^{++}\neq
\emptyset$.  This can also be viewed as a consequence of Gordan's
theorem (see Appendix A).  In general, $F$ is the convex hull of
the set of projections $\tau_i$ for which $Z(H)$ contains at least
one element of the form $a_1\tau_1+\cdots+a_n\tau_n$ with $a_i>0$.

When $H$ is finitely generated with basis $\{h_1,\dots,h_s\}$,
property $(B_m)$ for $H$ amounts to the solvability of a
system of Diophantine inequalities namely the conditions
that, for all $\epsilon>0$, all $\sigma\in\{-1,1\}$ and
all $i=1,\dots,s$, there exist integers $a_1,\dots,a_s\in\bZ$
such that all coordinates of $\sigma h_i-a_1h_1-\cdots-a_sh_s$
are at least $-\epsilon$. In Section \ref{sec:ex}, we give a
class of examples to which Theorem B applies, thereby solving
the corresponding system of Diophantine inequalities.

It follows from Theorem B that property $(B_m)$ is independent
of $m$.  In the next section, we extend this condition to an
arbitrary subgroup $H$ of a partially ordered abelian group $G$ with
order unit and we show that, in that context, it is again
independent of $m$.  From that point, we refer to it as property $(B)$.

In section \ref{sec:nec}, we establish a necessary condition
for $H$ to have property $(B)$ which, for $G=\bR^n$, reduces to
that of Theorem B.  We show in section \ref{sec:suf} that this
condition is also sufficient when the trace space of $G$ has
finite dimension, thereby completing the proof of Theorem B.
A construction in section \ref{sec:counterex} shows however
that this is not true for a general group $G$.  When $G$ is a simple
dimension group and $H$ is a convex subgroup of $G$ for which
$G/H$ is torsion-free, we prove in section \ref{sec:unperf}
that $G/H$ is unperforated if and only if $H$ has property $(B)$.
This complements \cite[Proposition B.1]{BeH} where the condition
is shown to be sufficient.

Section \ref{sec:ref}
answers a question of \cite{BeH} by giving necessary and
sufficient conditions for a trace of $G$ to be refinable when
$G$ is a simple dimension group with finitely many
pure traces.  The notion of refinable trace arose from a property
of measures on Cantor sets due to Akin \cite{Ak}, put in the
context of invariant probability measures on Cantor dynamical
systems, and subsequently translated to the setting of
dimension groups by S.~Bezuglyi and the first author
via the ordered $K_0$ functor \cite{BeH}.  Finally,
Appendix \ref{sec:A} explains the connection between Gordan's
theorem, Farkas's lemma, and some of our results.

\section{Condition $(B_m)$ in partially ordered abelian groups}
\label{sec:condBm}

By a \emph{partially ordered abelian group}, we mean an abelian
group $G$ equipped with a translation invariant partial order
$\le $.  The {\it positive cone} of such
a pair $(G,\le)$ is the set $G^+=\Set{x\in G}{0\le x}$.  It satisfies
\[
 (G^+)+(G^+)\subseteq G^+
 \quad\text{and}\quad
 (G^+)\cap(-G^+)=\{0\}.
\]
Conversely, any subset $G^+$ of an abelian group $G$ satisfying
these conditions makes $G$ into a partially ordered abelian group
upon defining, for $x,y\in G$, that $x\le y \ \Leftrightarrow\ y-x\in G^+$.
As usual, we write $x<y$ when $x\le y$ and $x\neq y$.

An \emph{order unit} (or \emph{strong unit}) of such a group $G$
is a nonzero element, $u$, of $G^{+}$ such that
for each $g \in G$, there exists a positive integer $N$
with $-Nu \leq g \leq Nu$.  The set of
all order units of $G$ is denoted $G^{++}$.

Let $(G,u)$ be a partially ordered group with order unit $u$.
A \emph{trace\/} (or {\it state\/}) on $G$ is a nonzero group
homomorphism $\tau\colon G\to \bR$ that is positive in the sense
that $\tau(G^+) \subseteq \bR^+$. We say that a trace $\tau$ is
\emph{normalized at $u$\/} if $\tau(u) = 1$, and we denote by
$S(G,u)$ the set of those traces.  It is a compact and convex
subset of $\bR^G$ for the product topology on $\bR^G$, see
[G; Proposition 6.2].  We denote by $\Aff S(G,u)$ the vector
space consisting of all convex-linear continuous real-valued
functions on $S(G,u)$, and we equip it with the supremum norm
so that it becomes a Banach space.  It is also
a partially ordered abelian group with respect to
pointwise ordering where, for $\varphi,\psi \in
\Aff S(G,u)$, we write $\varphi\le\psi$ when
$\varphi(\tau) \le \psi(\tau)$ for all $\tau\in S(G,u)$.

The \emph{affine representation\/} of $(G,u)$ is the second
dual map from $G$ to $\Aff S(G,u)$ which, to each $g\in G$,
associates the evaluation map $\widehat{g}\colon S(G,u)\to \bR$
given by $\widehat{g}(\tau) = \tau(g)$ for all $\tau\in S(G,u)$.
This is an order-preserving group homomorphism.  Within
$\Aff S(G,u)$, we identify $\bR$ with the subspace of constant
functions so that, for $g\in G$ and $a\in\bR$,
the condition $a\le\widehat{g}$ simply means $a\le\tau(g)$
for all $\tau\in S(G,u)$.  For much more on this, and
other aspects  of partially ordered abelian groups, the reader
is refered to \cite{G}.

We view $\bR^n$ as a partially ordered abelian group with respect to
the usual coordinatewise ordering.  Its positive cone is $(\bR^n)^+$
and its set of order units is $(\bR^n)^{++}$, as defined in the
introduction. In particular, the vector
\[
 \un=(1,\dots,1),
\]
with all coordinates equal to $1$, is an order unit of $\bR^n$.
The corresponding trace space $S(\bR^n,\un)$ is the simplex $K_n$ spanned by
the coordinate functions $\tau_1,\dots,\tau_n$ in $(\bR^n)^*$, and the
affine representation from $\bR^n$ to $\Aff K_n$ is an
isomorphism of vector spaces over $\bR$.  For $g\in \bR^n$
and $a\in\bR$, the condition $a\le\widehat{g}$
reduces to $a \le\tau_i(g)$ for each $i=1,\dots,n$.  Thus,
for a subgroup $H$ of $\bR^n$ and an integer $m\ge 2$,
condition $(B_m)$ can be restated as follows:
\begin{itemize}
\item[$(B_m)$\quad]
    for all $h \in H$ and $\epsilon > 0$, there exists
    $h' \in H$ such that $\widehat{h}-m\widehat{h'} \ge -\epsilon$.
\end{itemize}

We use this as a definition of the property
$(B_m)$ for a subgroup $H$ of $G$.  When it holds,
then so do ($B_n$) for all divisors $n>1$ of $m$,
and so do ($B_{m^j}$) for all integers $j\geq 1$.  More
generally, the next result shows that the conditions
$(B_m)$ with $m>1$ are mutually equivalent.

\begin{prop}
\label{condBm:prop:condB}
Let $(G,u)$ be as above.  Suppose that a subgroup $H$
of $G$ satisfies $(B_m)$ for some integer $m>1$.
Then $H$ satisfies $(B_n)$ for all integers $n > 1$.
\end{prop}

In view of this, we simply say that $H$ has property
$(B)$ if it satisfies $(B_m)$ for some (and thus
all) integers $m>1$.

\begin{proof} Let $h \in H$, let $n>1$ be an integer and let
$\epsilon>0$. Since $u$ is an order unit,
there exists $\ell\in\bN$ such that $h\le \ell u$.  Choose
an integer $j\ge 1$ such that $m^j\ge (2\ell+\epsilon)n/\epsilon$.
Since $H$ satisfies $(B_m)$, it satisfies ($B_{m^j}$), and so
there exists $h'\in H$ such that
\[
 \widehat{h}-m^j\widehat{h'}\ge -\epsilon/2.
\]
Since $h\le \ell u$, we have $\widehat{h}\le\ell$ and so the
above inequality yields
\[
 \widehat{h'}
 \le m^{-j}(\widehat{h}+\epsilon/2)
 \le m^{-j}(\ell+\epsilon/2)
 \le \epsilon/(2n).
\]
Writing $m^j=qn-r$ with integers $q\ge 1$ and $0\le r<n$,
we conclude that
\[
 \widehat{h}-nq\widehat{h'} \ge -r\widehat{h'}-\epsilon/2 \ge-\epsilon,
\]
thus $\widehat{h}-n\widehat{h''}\ge -\epsilon$ where $h''=qh'\in H$.
This shows that ($B_n$) is satisfied.
\end{proof}

If $G_0$ is a subgroup of $G$ containing $u$,
then $(G_0,u)$ is a partially ordered abelian group with
order unit, with positive cone $G_0^+=G^+\cap G_0$.
Importantly, by \cite[Theorem 3.2]{GH1} (see also
\cite[Corollary 4.3]{G}), the map $\rho\colon S(G,u)\to S(G_0,u)$
sending a trace on $G$ to its restriction to $G_0$ is a
{\it surjective} affine order-preserving homomorphism.
We deduce that the property $(B)$ for a subgroup $H$ of $G$
depends only on the induced ordering on $H+\bZ u$.

\begin{prop}
\label{condBm:prop:G0}
Let $(G,u)$ be as above, let $H$ be a subgroup of $G$, and
let $G_0=H+\bZ u$. Then $H$ has property $(B)$ within $(G,u)$
if and only if it has property $(B)$ within $(G_0,u)$.
\end{prop}

\begin{proof}
Let $m>1$ be an integer.  The condition $(B_m)$ for $H$ within $(G,u)$
requests that, for each $h\in H$ and each $\epsilon>0$, there
exists $h'\in H$ such that $\tau(h-mh')\ge -\epsilon$ for
all $\tau\in S(G,u)$.  The condition within $(G_0,u)$ is the same
except that $\tau$ varies in $S(G_0,u)$.  In view of the surjectivity
of the restriction map from $S(G,u)$ to $S(G_0,u)$, the two conditions
are thus the same.
\end{proof}

%
%

\section{A necessary condition for property $(B)$}
\label{sec:nec}

Let $H$ be a subgroup of a partially ordered abelian
group with order unit $(G,u)$.  The set
\[
 Z_G(H)
 := \Set{\tau \in S(G,u)}{\tau(h)= 0 \text{ for all $h \in H$}}
\]
is a compact convex subset of $S(G,u)$. When $G=\bR^n$ and
$u=\un$, this is the set denoted $Z(H)$ in Theorem B.

A \emph{face\/} of $S(G,u)$ is a (possibly empty) subset
$F$ of $S(G,u)$ such that any line segment in $S(G,u)$ whose
relative interior meets $F$, is contained in $F$.
For any subset $Z$ of $S(G,u)$, there is
a smallest face containing $Z$.
When $Z$ is convex (such as the set $Z_G(H)$ defined above),
it consists of all $\tau_1 \in S(G,u)$ for which there
exist $\tau_2 \in S(G,u)$ and $\lambda \in (0,1)$
such that $\lambda \tau_1 + (1-\lambda) \tau_2 \in Z$, see
[G, Proposition 5.7].  The \emph{extreme boundary} of $S(G,u)$,
denoted $\partial_e S(G,u)$, is the set of all traces
$\tau\in S(G,u)$, called \emph{pure traces}, which by themselves
constitute faces $\{\tau\}$ of $S(G,u)$.  For example,
$\partial_e S(\bR^n,\un)$ is the set of coordinates functions
$\{\tau_1,\dots,\tau_n\}$ in $(\bR^n)^*$.

The next result provides a necessary condition for property $(B)$
to hold. For $(G,u)=(\bR^n,\un)$, it is the same condition as in Theorem B.
The proof that we provide below is an adaptation of the arguments
from [BeH, Proposition B.4].

\begin{thm}
\label{nec:thm}
Let $(G,u)$ and $H$ be as above, and let $F$ be the
smallest face of $S(G,u)$ containing $Z_G(H)$.
Suppose that $H$ satisfies condition $(B)$. Then
for any $\tau$ in $F$, either $\tau(H)$ is zero, or
it is dense in $\bR$.
\end{thm}

\begin{proof} Pick $\tau_1$ in $F$. Since $Z_G(H)$ is a convex subset
of $S(G,u)$, there exist $\tau_2\in S(G,u)$ and
$\lambda \in (0,1)$ such that
$\lambda \tau_1 + (1-\lambda) \tau_2\in Z_G(H)$.
Then, for any $h$ in $H$, we have
\[
\tau_1(h) = - \theta \tau_2(h),
\quad \text{where $\theta = \frac{1-\lambda}{\lambda} > 0$}.
\]
Suppose that $\tau_1(H) = \bZ\delta$ for some real $\delta > 0$.
For every $h$ in $H$, it follows that
\[
 (\tau_1(h), \tau_2(h)) = \ell\delta (1,-\theta)
 \quad\text{for some $\ell\in\bZ$.}
\]
Setting $\epsilon = \delta \min\{1,\theta\}/2$,
we deduce that $(\tau_1(h), \tau_2(h)) \ge (-\epsilon,-\epsilon)$
in $\bR^2$ (for the componentwise ordering)
if and only if $\tau_1(h) = 0$.
Choose $h$ such that $\tau_1(h) = \delta$ generates $\tau_1(H)$,
and let $m>1$ be an integer.  Then, for every $h'$ in $H$,
we have $\tau_1(h) \neq m\tau_1(h')$, so $\tau_1(h-mh') \neq 0$
and by the above we obtain
\[
 (\widehat{h}-m\widehat{h'})(\tau_i)
  = \tau_i(h- mh')
  < -\epsilon
 \quad\text{for some $i\in\{1,2\}$}.
\]
This contradicts $(B_m)$. Hence $\tau_1(H)$ is either zero or dense in $\bR$.
\end{proof}

In section \ref{sec:suf}, we will show that the converse holds when
$S(G,u)$ spans a finite-dimensional subspace of $\bR^G$ or when $H$
has finite rank. This will complete the proof of Theorem B.  The next
section provides the last tool that we need for this purpose.

%
%

\section{An auxiliary result}
\label{sec:aux}

Throughout this section we fix an Euclidean space $E$ of finite
dimension $n>1$ with the scalar product of $\ux,\uy\in E$, denoted
$\ux\cdot\uy$.  We also fix a compact convex subset $K$ of $E$
containing $0$.  The notion of a face $F$ of $K$ and of the extreme
boundary $\partial_e K$ of $K$ is defined as in section \ref{sec:nec},
with $S(G,u)$ replaced by $K$.  In particular, there exists a
smallest face $F$ of $K$ containing $0$.  Our goal here
is to prove the following result.  Its relevance to
property $(B)$ will become clearer in the next section.

\begin{prop}
\label{aux:prop}
Let $F$ be the smallest face of $K$ containing $0$ and let
$Y$ be a subgroup of $E$ with $\bR Y=E$.  Suppose that
$\Set{\ux\cdot\uy}{\uy\in Y}$ is dense in $\bR$ for each
$\ux\in F\setminus\{0\}$.  Fix an arbitrary choice of\/
$\epsilon>0$ and of\/ $\uy\in Y$.  Then, there exists
$\uy_2\in Y$ such that $\ux\cdot(\uy-2\uy_2)\ge -\epsilon$
for each $\ux\in K$.
\end{prop}

For the rest of the section, we fix $F$ and $Y$ as in the
statement of the proposition. We also define
\[
 E''=\bR F  \et E'=(E'')^\perp,
\]
so that $E=E'\oplus E''$ is an orthogonal sum decomposition.
For the proof of the proposition, we will need the
following intermediate results.

\begin{lemma}
\label{aux:lemma:F}
The face $F$ is a neighbourhood of $0$ in $E''$.
\end{lemma}

\begin{proof}
This is clear if $F=\{0\}$ because then $E''=\{0\}$.
Suppose that $s=\dim_\bR E''$ is positive, and let
$\{\ux_1,\dots,\ux_s\}$ be a basis of $E''$ contained in $F$.
Since, for each $\ux\in F$ there exists $\lambda>0$ such
that $-\lambda\ux\in F$, we may assume that $-\ux_1,\dots,
-\ux_s\in F$.  Then $F$ contains the convex hull of the
$2s$ points $\pm\ux_1,\dots,\pm\ux_s$ which is a
neighbourhood of $0$ in $E''$.
\end{proof}

\begin{lemma}
\label{aux:lemma:K'}
Let $K'=\proj_{E'}K$ denote the image of $K$ under
the orthogonal projection on $E'$.  Then $K'$ is a compact
convex subset of $E'$ with $0\in\partial_eK'$.
\end{lemma}

\begin{proof}
Since the orthogonal projection on $E'$ is linear, thus
continuous, the image $K'$ of $K$ is convex, compact, and
contains $0$.  If $0\notin\partial_eK'$, there exists
$\ux'\in K'\setminus\{0\}$ such that $-\ux'\in K'$.  We can
write $\ux'=\ux_1+\uy_1$ and $-\ux'=\ux_2+\uy_2$ for some
$\ux_1,\ux_2\in K$ and some $\uy_1,\uy_2\in E''$.  By Lemma 3, there
exists $\delta\in(0,1/2)$ such that $2\delta\uy_i \in
F \subseteq K$ for $i=1,2$.  Since $K$ is convex,
containing $0$ and $\ux_i$, it also contains $2\delta\ux_i$
for $i=1,2$. So it contains $\delta\ux'$ and $-\delta\ux'$,
which in turn implies that $\delta\ux'\in F$.  However, this
is impossible since $F\cap E'\subseteq E''\cap E'=\{0\}$.
This contradiction shows that $0\in\partial_eK'$.
\end{proof}

\begin{lemma}
\label{aux:lemma:Y}
The orthogonal projection $\proj_{E''}Y$ of $Y$ on $E''$
is dense in $E''$.
\end{lemma}

\begin{proof}
The group $\proj_{E''}Y$ is dense in $E''$ if and only if
the orthogonal projection of $Y$ on $\bR\ux$ is dense in
$\bR\ux$ for each $\ux\in E''\setminus\{0\}$ or,
equivalently, if and only if $\Set{\ux\cdot\uy}{\uy\in Y}$
is dense in $\bR$ for each $\ux\in E''\setminus\{0\}$.
Since, by Lemma \ref{aux:lemma:F}, $F$ is a neighbourhood
of $0$ in $E''$, this is equivalent to the hypothesis that
$\Set{\ux\cdot\uy}{\uy\in Y}$ is dense in $\bR$  for each
$\ux\in F\setminus\{0\}$.
\end{proof}

\medskip
The next lemma is a basic tool of inhomogeneous Diophantine
approximation.

\begin{lemma}
\label{aux:lemma:convex}
Let $C$ be a convex neighbourhood of $0$
in a real vector space $V$ of finite dimension $s\ge 1$,
and let $\{\uv_1,\dots,\uv_s\}$ be a basis of $V$
contained in $C$.  Then any translate of $sC$ in $V$
contains at least one element of the group
$\bZ\uv_1+\cdots+\bZ\uv_s$.
\end{lemma}

\begin{proof}
Let $\uv\in V$.  Write $\uv=a_1\uv_1+\cdots+a_s\uv_s$
with $a_1,\dots,a_s\in\bR$.  Then $\uv+sC$ contains the
point $\lceil a_1\rceil\uv_1+\cdots+\lceil a_s\rceil\uv_s$,
where $\lceil a\rceil$ stands for the least integer greater
than or equal to $a$.
\end{proof}

It has the following consequence.

\begin{lemma}
\label{aux:lemma:neighbourhood}
For each neighbourhood $U''$ of\/ $0$ in $E''$, there exists
a compact neighbourhood $U'$ of\/ $0$ in $E'$ such that
each translate of\/ $U'+U''$ in $E$ contains at least one element
of $Y$.
\end{lemma}

\begin{proof}
It suffices to prove this for a convex
neighbourhood $U''$ of $0$, assuming that $s=\dim_\bR(E'')>0$.
Since, by Lemma \ref{aux:lemma:Y},
$\proj_{E''}Y$ is dense in $E''$, there exist elements
$\uy_1,\dots,\uy_s$ of $Y$ whose projections on $E''$ form a
basis of $E''$ contained in $(ns)^{-1}U''$.  Let
$Y''=\bZ\uy_1+\dots+\bZ\uy_s$.  By Lemma \ref{aux:lemma:convex},
each translate of $n^{-1}U''$ contains at least one element
of $\proj_{E''}Y''$.
In particular, each element of $Y$ is congruent modulo $Y''$
to a point $\uy$ with $\proj_{E''}\uy\in n^{-1}U''$.  Since
$\bR Y=E$, we may therefore complete $\{\uy_1,\dots,\uy_s\}$
to a basis $\{\uy_1,\dots,\uy_n\}$ of $E$ contained in $Y$,
with $\proj_{E''}\uy_j\in n^{-1}U''$ for $j=1,\dots,n$.
Choose a compact neighbourhood $U'$ of $0$ in $E'$ such that $n^{-1}U'$
contains the projections of $\uy_1,\dots,\uy_n$ on $E'$.
Then $n^{-1}(U'+U'')$ contains $\uy_1,\dots,\uy_n$ and so,
by Lemma \ref{aux:lemma:convex}, each translate of $U'+U''$
in $E$ contains at least one element of $\bZ\uy_1+\cdots+\bZ\uy_n$.
\end{proof}

\begin{lemma}
\label{aux:lemma:ball}
Let $\epsilon>0$ and $\rho>0$ be arbitrary. For each
nonzero subspace $V$ of $E'$, there exists
a closed ball $B$ of\/ $V$ of radius $\rho$ such
that $\ux\cdot\uy\ge -\epsilon$ for all $\ux\in K'\cap V$
and all $\uy\in B$.
\end{lemma}

\begin{proof}
Let $V$ be a nonzero subspace of $E'$.
By Lemma \ref{aux:lemma:K'}, the set $K'$ is a compact
convex subset of $E'$ with $0\in \partial_e K'$.
Then, $K'\cap V$ is a compact convex subset of $V$ with
$0\in\partial_e (K'\cap V)$.  In particular, $0$ belongs
to the topological boundary of $K'\cap V$ in $V$.
So, there exists a unit vector $\uu$ of $V$ such that
$\uu\cdot\ux\ge 0$ for each $\ux\in K'\cap V$,
by \cite[Proposition 5.10]{G}.  We prove
the existence of $B$ by induction on $s=\dim_\bR V$.

If $s=1$, we have $K'\cap V\subseteq \bR^+\uu$.  Then $B=[0,2\rho]\uu$
is a closed ball of $V$ of radius $\rho$ with $\ux\cdot\uy\ge 0$
for all $\ux\in K'\cap V$ and $\uy\in B$.

Suppose now that $s>1$.  Define
$V_0=\Set{\ux\in V}{\uu\cdot\ux=0}$ and $K_0=K'\cap V_0$.
Then, $V_0$ is a subspace of $V$ of dimension $s-1$.
So, we may assume the existence of a
closed ball $B_0$ of $V_0$ of radius $\rho$ such that
$\ux\cdot\uy\ge -\epsilon/2$ for all $\ux\in K_0$
and all $\uy\in B_0$.  Put
\[
 L=\sup\Set{\|\ux\|}{\ux\in K'\cap V}
 \quad
 \text{and}
 \quad
 M=\sup\Set{\|\uy\|}{\uy\in B_0}.
\]
Define also $U=\Set{\ux\in V}{\|\ux\|<\epsilon/(2M)}$ and
$F_\delta=\Set{\ux\in K'\cap V}{\uu\cdot\ux\le \delta}$ for each $\delta>0$.
Then $\Set{F_\delta}{\delta>0}$ is a collection of closed subsets
of $K'\cap V$, stable under finite intersections, whose intersection is $K_0$.
Since $K'\cap V$ is compact and since $K_0+U$ is an open subset of $V$
containing $K_0$, there exists therefore $\delta>0$ for which
$F_\delta\subseteq K_0+U$.  Let $B$ be any ball of $V$ of radius $\rho$
contained in $B_0+[R,\infty)\uu$, where $R=ML/\delta$.  We claim that
$B$ has the required property.

To show this, choose any $\ux\in K'\cap V$ and $\uy\in B$.  We may write
$\uy=\uy_0+t\uu$ where $\uy_0\in B_0$ and $t\ge R$.  If $\ux\notin
F_\delta$, we have $\uu\cdot\ux>\delta$ and so
\[
 \ux\cdot\uy
  = \ux\cdot\uy_0+t\uu\cdot\ux
  > -LM+R\delta =0.
\]
Otherwise, we have $\ux\in F_\delta\subseteq K_0+U$, so there
exists $\ux_0\in K_0$ such that $\|\ux-\ux_0\|<\epsilon/(2M)$.
Since $\uu\cdot\ux\ge 0$ and $\|\uy_0\|\le M$, we obtain
\[
 \ux\cdot\uy
   =\ux\cdot\uy_0+t\uu\cdot\ux
   \ge \ux\cdot\uy_0
   =\ux_0\cdot\uy_0+(\ux-\ux_0)\cdot\uy_0
   \ge -\frac{\epsilon}{2}-\frac{\epsilon}{2M}M
   = - \epsilon.
\]
\end{proof}

\begin{proof}[Proof of Proposition \ref{aux:prop}]
By Lemma \ref{aux:lemma:Y}, the group $\proj_{E''}Y$ is
dense in $E''$.  Thus, the same is true of $\proj_{E''}(2Y)$.
If $E'=\{0\}$, this means that $2Y$ is dense in $E$ and
the conclusion follows.  We thus assume that $E'\neq\{0\}$.
Then, for each $\delta>0$, Lemma \ref{aux:lemma:neighbourhood}
shows the existence of $\rho>0$ such that each translate
of $B'_\rho+B''_\delta$ contains at least one element of $2Y$,
where $B'_\rho=\Set{\uy'\in E'}{\|\uy'\|\le \rho}$ and
$B''_\delta=\Set{\uy''\in E''}{\|\uy''\|\le \delta}$.
Moreover, applying Lemma \ref{aux:lemma:ball} to the choice
of $V=E'$, we obtain a translate $B'$ of $B'_\rho$
with the property that $\ux'\cdot\uy'\ge -\epsilon/2$ for any
$\ux'\in K'$ and any $\uy'\in B'$.

We apply these observations with $\delta=\epsilon/(2L)$ where
$L=\sup\{\|\ux\|\,;\, \ux\in K\}$.  Put $B=B'+B''_\delta$ for the
corresponding choice of $B'$.  Then any translate of $B$ in $E$
contains at least one element of $2Y$.  In particular, for the
given $\uy\in Y$, there exists $\uy_2\in Y$ such that
$2(-\uy_2)\in(-\uy)+B$, and thus $\uy-2\uy_2\in B$.  Write
$\uy-2\uy_2=\uy'+\uy''$ with $\uy'\in B'$ and $\uy''\in B''_\delta$.
Then, for any $\ux\in K$, we find
\[
\ux\cdot(\uy-2\uy_2)
 = \proj_{E'}(\ux)\cdot\uy'+\proj_{E''}(\ux)\cdot\uy''
 \ge (-\epsilon/2)-L\delta \ge -\epsilon
\]
as requested.
\end{proof}

%
%

\section{Sufficient conditions for property (B)}
\label{sec:suf}

We start with the following observation.

\begin{prop}
\label{suf:prop:ZHempty}
Let $H$ be a subgroup of a partially ordered abelian group with
order unit $(G,u)$.  We have
\[
 Z_G(H)=\emptyset
 \quad \Longleftrightarrow \quad
 H\cap G^{++}\neq \emptyset.
\]
When this happens, the group $H$ has property $(B)$.
\end{prop}

The first assertion generalizes Gordan's theorem \cite{Gor},
see Appendix \ref{sec:A}.  When $G=\bR^n$ and $u=\un$, the
second proves Theorem B in the case that $Z_G(H)=Z(H)$ is empty.

\begin{proof}
If $H$ contains no order
units, then $G_0=H\oplus\bZ u$ is a direct sum and the map
$\tau_0\colon G_0\to \bR$ given by $\tau_0(h+\ell u)=\ell$
for each pair $(h,\ell)\in H\times\bZ$ is a positive homomorphism
for the restriction to $G_0$ of the partial order on $G$.
Since $u\in G_0$ and $\tau_0(u)=1$, this maps extends to a trace
$\tau\in S(G,u)$.  Then $Z_G(H)$ is not empty, as it contains $\tau$.

Conversely, if $H$ contains an order unit $v$, then
$u\le k v$ for some positive integer $k$,
so $\widehat{v} \ge k^{-1}\widehat{u} = k^{-1}$,
and thus $Z_G(H)=\emptyset$.  Moreover, let $m>1$ be an integer
and let $h\in H$.  Since $mv$ is an order unit, we have
$-h\le \ell mv$ for some $\ell\in\bN$, thus $\widehat{h}-m\widehat{h'}\ge 0$
for $h'=-\ell v\in H$. Thus $H$ satisfies $(B_m)$ as well.
\end{proof}

\begin{thm}
\label{suf:thm}
Let $(G,u)$ be a partially ordered abelian group with order unit,
let $H$ be a subgroup of $G$, and let $F$ be the smallest face
of $S(G,u)$ containing $Z_G(H)$.
Let $\widehat{H}$ denote the canonical image of $H$ in $\Aff(S(G,u))$.
Suppose that $\bR\widehat{H}$ is finite-dimensional
and that $\tau(H)$ is dense in $\bR$ for each $\tau\in
F\setminus Z_G(H)$. Then, $H$ has property $(B)$.
\end{thm}

In particular, $H$ has property $(B)$ if $Z_G(H)$ is a face
of $S(G,u)$ and $\bR\widehat{H}$ has finite dimension.  We will
show in the next section that property $(B)$ may fail if
finite-dimensionality of $\bR\widehat{H}$ is dropped.

\begin{proof}
We may assume that $Z_G(H)\neq\emptyset$ since otherwise we already
know, by Proposition \ref{suf:prop:ZHempty}, that $H$ has property $(B)$.
From this, we proceed in two steps.  We first prove the statement
when $G=H+\bZ u$.  Then, we show that the general case
reduces to this special case.

So, assume first that $G=H+\bZ u$.  Then, as $Z_G(H)$ is not empty,
we have $H\cap\bZ u=\{0\}$ and $Z_G(H)$ consists of the single
group homomorphism $\tau_0\colon G \to\bZ$ given by $\tau_0(h+nu)=n$
for each $h\in H$ and each $n\in\bZ$.
Choose $h_1,\dots,h_n\in H$ such that $\{\widehat{h}_1,\dots,\widehat{h}_n\}$
is a basis of $\bR\widehat{H}$ over $\bR$.  We denote by
$\varphi\colon H\to \bR^n$ the group homomorphism that sends
an element of $h$ to the coordinates of $\widehat{h}$ in that basis, namely
the $n$-tuple $\varphi(h)=(y_1,\dots,y_n)\in\bR^n$ such that
\[
 \widehat{h}=y_1\widehat{h}_1+\cdots+y_n\widehat{h}_n.
\]
We also form the linear map $\psi\colon \bR^G\to \bR^n$ given by
\[
 \psi(\tau)=(\tau(h_1),\dots,\tau(h_n))
\]
for each $\tau\in R^G$.  Then, for each $\tau\in S(G,u)$ and
each $h\in H$, we have
\[
 \tau(h)=\widehat{h}(\tau)=\psi(\tau)\cdot\varphi(h)
\]
using the standard scalar product on $\bR^n$.  We deduce from this that
$\psi$ is one-to-one on $S(G,u)$, because if elements $\tau_1,\tau_2$ of $S(G,u)$
have the same image under $\psi$ then, by the formula above, they
coincide on $H$ and therefore they coincide on $G$ (since they
take the value $1$ at $u$).

Let $K=\psi(S(G,u))$ and $Y=\varphi(H)$.  Then,
$K$ is a convex subset of $\bR^n$ containing $\psi(\tau_0)=0$ and $\psi$
induces an affine homeomorphism from $S(G,u)$ to $K$.  Accordingly,
$F_0:=\psi(F)$ is the smallest face of $K$ containing $0$.  Moreover,
$Y$ is a subgroup of $\bR^n$ which contains $\bZ^n$, the image of
$\bZ h_1+\cdots+\bZ h_n$ under $\varphi$.  So, we have $\bR Y=\bR^n$.
The hypothesis also implies that $\{\ux\cdot\uy\,;\,\uy\in Y\}$ is dense
in $\bR$ for each $\ux\in F_0\setminus \{0\}$.  So, $K$ and $Y$ fulfil
all the hypotheses of Proposition 2 within the Euclidean space $\bR^n$.

Let $\epsilon>0$ and $h\in H$. Putting $\uy=\varphi(h)$, Proposition 2
yields a point $\uy'\in Y$ such that $\ux\cdot(\uy-2\uy')\ge -\epsilon$
for each $\ux\in K$. Choose any $h'\in H$ such that $\varphi(h')=\uy'$.
Then, the above property translates into $\widehat{h}(\tau)-2\widehat{h'}(\tau)\ge
-\epsilon$ for each $\tau\in S(G,u)$, showing that $H$ has property
$(B_2)$.

Now, we turn to the general case.  Consider the group $G_0=H+\bZ u$
with the induced ordering from $G$. Since the map
$\rho\colon S(G,u)\to S(G_0,u)$ sending a trace on $(G,u)$ to
its restriction on $(G_0,u)$ is a surjective affine
order-preserving homomorphism (see section \ref{sec:condBm}),
we have
\[
 Z_G(H)=\rho^{-1}(Z_{G_0}(H)) \quad\text{and}\quad F=\rho^{-1}(F_0),
\]
where $F_0$ denotes the smallest face of $S(G_0,u)$ containing $Z_{G_0}(H)$.
In particular $Z_{G_0}(H)$ is not empty and the hypothesis that $\tau(H)$
be dense in $\bR$ for each $\tau\in F\setminus Z_G(H)$ implies
the same for each $\tau\in F_0\setminus Z_{G_0}(H)$.
Let $\widehat{H}_0$ denote the canonical image of $H$ in $\Aff(S(G_0,u))$.
The dual map $\rho^*$ from $\Aff(S(G_0,u))$ to $\Aff(S(G,u))$ given by
composition with $\rho$ is $\bR$-linear and restricts to an
isomorphism of $\bR$-vector spaces from $\bR\widehat{H}_0$ to $\bR\widehat{H}$.
In particular, $\bR\widehat{H}_0$ is finite-dimensional.  Thus the hypotheses of
the theorem hold for $H$ as a subgroup of $G_0$ and therefore,
by the special case proved above, the group $H$ has property $(B)$
within $G_0$.  By Proposition \ref{condBm:prop:condB}, we conclude that
it also has property $(B)$ within $G$.
\end{proof}

%
%

\section{A class of counter-examples}
\label{sec:counterex}

The criterion for property $(B)$ given in Theorem \ref{suf:thm} requires
that $\bR\widehat{H}$ have finite dimension.  The goal
of this section is to exhibit examples showing that the
theorem is false without this hypothesis.
Our construction is based on the following observation.

\begin{prop}
\label{counterex:prop1}
Let $X$ be any infinite set, and let $H$ be the subgroup
of $\bR^X$ generated by an infinite sequence of bounded
functions $(h_i\colon X\to\bR)_{i\in\bN}$ which is linearly
independent over $\bR$.  Suppose that, for each
nonzero $f$ in $\bR H$, there exists $x\in X$ for which
$f(x)>0$.  Then, there is a sequence of
positive integers $(m_i)_{i\in\bN}$ such that each
nonzero element $h$ of\/ $H'=\sum_{i\ge 1}\bZ m_ih_i$ satisfies
$\sup_X h \ge 1$.
\end{prop}

\begin{proof}
We construct recursively the integers $m_1,m_2,\dots$ so that
$\sup_X h \ge 1$ for each $n\ge 1$ and each nonzero $h$ in
$H_n = \bZ m_1h_1+\cdots+\bZ m_nh_n$.

For $n=1$, we choose $x_1,x_2\in X$ such that $h_1(x_1)>0$
and $(-h_1)(x_2)>0$, and we select $m_1\in\bN$ large enough
so that $m_1h_1(x_1)\ge 1$ and $-m_1h_1(x_2)\ge 1$.  Then, for each
nonzero $h\in H_1=\bZ m_1h_1$, we obtain
$\sup_X h \ge \max\{h(x_1),h(x_2)\} \ge 1$.

Suppose $m_1,\dots,m_n$ constructed for some integer $n\ge 1$.
For each point $\uu=(u_1,\dots,u_{n+1})$ in the unit sphere $S$
of $\bR^{n+1}$, the function $f=u_1h_1+\cdots+u_{n+1}h_{n+1}$
is a nonzero element of $\bR H$, and so there exists
$x_\uu\in X$ such that $f(x_\uu)>0$.  Define
\[
 V_\uu
 =
 \{(v_1,\dots,v_{n+1})\in\bR^{n+1}\,;\,
   (v_1h_1+\cdots+v_{n+1}h_{n+1})(x_\uu)>0\}.
\]
Then $(V_\uu)_{\uu\in S}$ is an open covering of $S$.  Since
$S$ is compact, it admits a finite subcover by open sets
corresponding to points $x_1,\dots,x_k\in X$.  This means
that the function $g\colon S\to \bR$ given by
\[
 g(u_1,\dots,u_{n+1})
 =
 \max_{1\le i\le k}
 (u_1h_1+\cdots+u_{n+1}h_{n+1})(x_i)
\]
is strictly positive on $S$.  Since it is continuous, it is
therefore bounded below by some positive constant $\delta>0$.
We chose $m_{n+1}\in\bN$ so that $m_{n+1}\delta\ge 1$, and
claim that this integer has the desired property.

To show this, choose any nonzero element $h$ of $H_{n+1}=
H_n+\bZ m_{n+1}h_{n+1}$.  If $h\in H_n$, then by hypothesis we
have $\sup_X h \ge 1$.  Otherwise, we can write
$h=\ell_1m_1h_1+\cdots+\ell_{n+1}m_{n+1}h_{n+1}$ for some
integers $\ell_1,\dots,\ell_{n+1}\in\bZ$ with $\ell_{n+1}\neq 0$.
Put $\uv=(\ell_1m_1,\dots,\ell_{n+1}m_{n+1})$ and
$\uu=\|\uv\|^{-1}\uv$. Since $\uu\in S$ and $\|\uv\|\ge m_{n+1}$,
we obtain
\[
 \sup_X h
   \ge \max_{1\le i\le k} h(x_i)
    = \|\uv\| g(\uu)
   \ge m_{n+1} \delta
   \ge 1,
\]
as desired.
\end{proof}

Consider the vector space $G=\cC(X,\bR)$ consisting of continuous
real-valued functions on a compact Hausdorff topological space
$X$.  This is a partially ordered abelian group with respect to the
pointwise ordering and, as such, it admits the constant
function $\un$ as an order unit.  Each trace $\tau$ in
$S(G,\un)$ is given by $\tau(f)=\int_X f\,d\mu$ for a unique
Borel probability measure $\mu$ on $X$ and this identifies
$S(G,\un)$ to the set $M^+_1(X)$ of probability measures
on $X$ with the weak* topology (see \cite[Chapter 5]{G}).
Moreover, the extreme boundary of $S(G,\un)$ is homeomorphic
to $X$ under the map that sends a point $x$ in $X$ to
the evaluation $\epsilon_x$ at $x$ given by
\[
 \epsilon_x(f)=f(x)
 \quad
 (f\in G)
\]
corresponding to the point-mass measure $\delta_x$
at $x$ \cite[Proposition 5.24]{G}.  In particular $\{\epsilon_x\}$
is a face of $S(G,\un)$ for any $x\in X$.

\begin{prop}
Let $X$ be an infinite compact metrizable topological space.
Choose $x_0\in X$ such that
$X_0:=X\setminus\{x_0\}$ is dense in $X$, and consider $G=\cC(X,\bR)$
as a partially ordered abelian group as above.  Then
there exists a subgroup $H$ of $G$ with the following
properties:
\begin{itemize}
\item[(i)] $Z_G(H)=\{\epsilon_{x_0}\}$, \par
\item[(ii)] $\sup_X h \ge 1$ for each $h\in H\setminus\{0\}$.
\end{itemize}
Then $Z_G(H)$ is a face of $S(G,\un)$, and $H$ does not satisfy
property $(B)$.
\end{prop}

The existence of $x_0$ follows from the
fact that $X$ is not discrete, since a discrete compact Hausdorff
space is finite.

\begin{proof}
Let $G_0=\ker(\epsilon_{x_0})=\{g\in G\,;\, g(x_0)=0\}$.  Since
$X$ is compact metrizable, $\cC(X,\bR)$ is a separable
topologial space; so the group $G_0$
contains a dense countable sequence $(f_n)_{n\in\bN}$ consisting
of continuous functions with compact support in $X_0$.  Choose
a dense sequence $(x_n)_{n\in\bN}$ in $X_0$ and a subsequence
$(y_n)_{n\in\bN}$ converging to $x_0$.  Let
\[
 \mu = \sum_{n=1}^\infty \frac{1}{2^n}\delta_{x_n}
       + \sum_{n=1}^\infty \delta_{y_n}.
\]
This defines
a finite positive measure on each compact subset $K$ of $X_0$ since
$X\setminus K$ contains $y_n$ for all but finitely many $n$.
Since $\mu(X_0)=\infty$, there is also, for each
$m\in\bN$, a function $g_m\in G_0$ with compact support
in $X_0$ such that $\|g_m\|_\infty\le 1/m$ and $\int_X g_m d\mu=1$.
Let $H_0$ be the subgroup of $G_0$ generated by the functions
\[
 f_n - \left(\int_X f_n\, d\mu \right)g_m,
 \qquad (m,n)\in\bN^2.
\]
Since the sequence $(g_m)_{m\in\bN}$ converges uniformly to
$0$ on $X$, the topological closure $\overline{H}_0$ of $H_0$
in $G$ contains $f_n$ for each $n\in\bN$.  We conclude that
$\overline{H}_0=G_0$, and so $Z_G(H_0)=Z_G(G_0)=\{\epsilon_{x_0}\}$.
Moreover, we have $\int_X f\, d\mu=0$ for each $f$ in $H_0$.
Therefore, the same is true for each $f$ in $\bR H_0$ and thus,
for each non-zero $f\in\bR H_0$, there exists $n\in\bN$ for
which $f(x_n)>0$.  Choose a basis $(h_i)_{i\in\bN}$ of
$\bR H_0$ contained in $H_0$.  By Proposition \ref{counterex:prop1}, there exists
a sequence $(m_i)_{i\in\bN}$ in $\bN$ such that the subgroup
$H$ of $G$ spanned by $(m_ih_i)_{i\in\bN}$ satisfies condition
(ii).  It also satisfies condition (i) since $Z_G(H)=Z_G(H_0)$.

Finally, let $h\in H\setminus 2H$ (for example $h=m_1h_1$),
and let $h'\in H$.  Then $2h'-h\neq 0$ satisfies
$\sup_X(2h'-h)\ge 1$, so $h(x)-2h'(x)\le -1$ for some $x\in X$,
meaning that $\tau(h)-2\tau(h')\le -1$ for the trace $\tau=\epsilon_x$.
Thus $H$ does not satisfy property $(B_2)$.
\end{proof}

%
%

\section{Link with unperforation of quotients}
\label{sec:unperf}

Let $G$ be a partially ordered abelian group.  We recall the following
definitions:
\begin{itemize}
\item $G$ is \emph{directed} if any finite subset
  of $G$ has an upper bound in $G$;
\item $G$ is \emph{simple} if it is nonzero, directed,
  and $G^+\setminus\{0\}=G^{++}$;
\item $G$ is \emph{unperforated} if the condition $mg\ge0$
  with $g\in G$ and $m\in\bN$ implies that $g\ge 0$;
\item $G$ has the \emph{Riesz interpolation property} if,
  given $g_1,g_2,g_1',g_2'\in G$ with $g_i\le g'_j$
  for each $i,j\in\{1,2\}$, there exists $g\in G$ such
  that $g_i\le g\le g_j'$ for each $i,j\in\{1,2\}$.
\end{itemize}
Clearly, $G$ is directed if it admits an order unit.
If $G$ is unperforated, then it is torsion-free as an
abelian group, and an element $g$ of $G$ is an order
unit if and only if $\tau(g)>0$ for all traces $\tau$ on $G$
\cite[Theorem 1.4]{EHS} or \cite[Corollary 4.13]{G}.
The group $G$ is called a \emph{dimension group} if it is
directed, unperforated, and has the Riesz interpolation
property.  For example, $\bR^n$ is a simple
dimension group for the \emph{strict ordering} with
positive cone $\{0\}\cup (\bR^n)^{++}$, and so is any dense
subgroup $G$ of $\bR^n$ for the inherited ordering.
This also applies to $\bZ$ but not to $\bZ^n$ if $n>1$.

A subgroup $H$ of $G$ is said to be \emph{convex} if whenever
$h\le g\le h'$ with $h,h'\in H$ and $g\in G$, we have $g\in H$.
Suppose that $H$ is such a subgroup.
Then the quotient $G/H$ is a partially ordered abelian group
with positive cone $(G/H)^{+}=(G^{+}+H)/H$\,; given $g,g'\in G$,
we have $g+H\le g'+H$ if and only if $g\le g'+h$ for some $h\in H$.
The following result links
unperforation of $G/H$ to property $(B)$ for $H$.

\begin{prop}
\label{unperf:prop}
Suppose that $(G,u)$ is a simple unperforated partially ordered
abelian group with order unit, and let $H$ be a convex subgroup
of $G$ for which $G/H$ is torsion-free.
\begin{itemize}
\item[(i)] If $H$ has property $(B)$, then $G/H$ is unperforated.
\item[(ii)] If $G/H$ is unperforated and if the set
  $\Set{\widehat{g}}{g\in G}$ is dense in $\Aff S(G,u)$, then
  $H$ has property $(B)$.
\end{itemize}
\end{prop}

The argument for (i) follows the proof of [BeH, Proposition B.1]
together with a simplification.

\begin{proof} (i) Suppose that $H$ has property $(B)$.
Let $g\in G$ and $m\ge 2$ be an integer such that $m(g+H)\ge H$.
We need to show that $g+H\ge H$. If $m(g+H)=H$ then
$g+H=H$ because $G/H$ is torsion-free, and we are done.
Otherwise, we have $mg+h>0$ for some $h\in H$ and so
$mg+h\in G^{++}$ because $G$ is simple.  Choose $n\in\bN$
such that $n(mg+h)\ge u$.  Then, we have
$m\widehat{g}+\widehat{h} \ge 1/n$.  Moreover,
since $H$ satisfies $(B_m)$, there exists $h'\in H$
with $-\widehat{h}+m\widehat{h}'\ge -1/(2n)$.  Combining
these inequalities, we deduce that $\widehat{g}+\widehat{h}'
\ge 1/(2mn)$.  This implies that $g+h'\in G^{++}$ since $G$
is unperforated, and so $g+H\ge H$.

(ii) Suppose that $G/H$ is unperforated and that
$\Set{\widehat{g}}{g\in G}$ is dense in $\Aff S(G,u)$.
Let $h\in H$ and $\epsilon>0$.  The set of $a\in \Aff S(G,u)$
with $\epsilon/4 < a-\widehat{h}/2 <\epsilon/2$ is
open and not empty as it contains $\widehat{h}/2+\epsilon/3$.
Thus it contains $\widehat{g}$ for some $g\in G$. This means
that $\epsilon/2 \le 2\widehat{g}-\widehat{h}\le \epsilon$.
The first inequality
implies that $2g-h\in G^{++}$ since $G$ is unperforated.
In particular, we have $2g-h\ge 0$, so $2(g+H)\ge H$, and
unperforation of $G/H$ yields $g + H \in (G/H)^+$.
Thus there exists $h'$ in $H$ such that $g-h' \in G^+$.
Then $\widehat{g}\ge\widehat{h}'$ and we obtain $-\epsilon\le
\widehat{h}-2\widehat{g}\le \widehat{h}-2\widehat{h}'$,
showing that $H$ satisfies $(B_2)$.
\end{proof}

\begin{cor}
\label{unperf:cor}
Suppose that $(G,u)$ is a simple dimension group with order unit,
and let $H$ be a convex subgroup of $G$ for which $G/H$ is torsion-free.
Then $G/H$ is unperforated if and only if $H$ has property $(B)$ inside
$G$.
\end{cor}

Necessity of the condition is proved in [BeH, Proposition B.1].

\begin{proof}
If $G$ is cyclic as an abelian group, then $H$ is $\{0\}$ or $G$.  So
$G/H$ is unperforated and $H$ has property $(B)$. If $G$ is noncyclic, then
$\Set{\widehat{g}}{g\in G}$ is dense in $\Aff S(G,u)$ by
\cite[Corollary 4.10]{GH2} or \cite[Theorem 14.14]{G}, and the conclusion
follows from the Proposition.
\end{proof}

An example of a convex subgroup $H$ of a simple dimension group $G$ such that
$G/H$ is torsion-free but \emph{holey}, i.e., not unperforated, is given in
\cite[Appendix B]{BeH}.  It suffices to take for $G$ a rank $3$ dense
subgroup of $\bR^2$ of the form $G=\bZ^2+\bZ(\alpha,\beta)$ with the strict
ordering, where $\{1, \alpha,\beta \}$ is linearly independent over $\bQ$,
and to take $H=\bZ(-1,1)$ (it is convex since $H\cap G^+=\{(0,0)\}$).

%
%

\section{Application to refinable traces}
\label{sec:ref}

The motivating reason to consider unperforation of quotients comes from the
study of refinable traces (originally refinable measures on Cantor
dynamical systems, \cite{Ak}). We first recall some definitions from
\cite[Sections 1 and 7]{BeH}.

Let $(G,u)$ be a dimension group with order unit, and let $U$ be a
nonempty subset of $S(G,u)$.  We say that $U$ is
\begin{itemize}
\item \emph{refinable} if, whenever $a_1, a_2, b\in G^+$ satisfy
  $\tau(a_1)+\tau(a_2)=\tau(b)$ for each $\tau\in U$, there exist
  $a_1', a_2' \in G^+$ such that $a_1'+a_2'=b$ and
  $\tau(a_i')=\tau(a_i)$ for each $i=1,2$ and each $\tau\in U$;
\item \emph{weakly good} if, for any $a,b\in G^+\setminus\{0\}$
  satisfying
  \[
   \inf_{\tau\in U} (\tau(b)-\tau(a))>0
  \]
  there exists $a'\in G^+\setminus\{0\}$ such that $a'<b$ and
  $\tau(a')=\tau(a)$ for each $\tau\in U$;
\item \emph{good} if, for any $a\in G$ and $b\in G^+\setminus\{0\}$
  satisfying
  \[
   \inf_{\tau\in U} \tau(a) > 0
   \quad\text{and}\quad
   \inf_{\tau\in U} (\tau(b)-\tau(a)) > 0
  \]
  there exists $a'\in G^+\setminus\{0\}$ such that $a'<b$ and
  $\tau(a')=\tau(a)$ for each $\tau\in U$.
\end{itemize}
Note that, when $U$ is a compact subset of
$S(G,u)$, for example when $U$ is finite or when $U=Z_G(H)$
for a subgroup $H$ of $G$, the condition $\inf_{\tau\in U}
\tau(a) > 0$ for $a\in G$ is equivalent to asking that
$\tau(a)>0$ for each $\tau\in U$; similarly the
condition $\inf_{\tau\in U} (\tau(b)-\tau(a)) > 0$
for $a,b\in G$ is then equivalent to $\tau(a)<\tau(b)$ for
each $\tau\in U$.

We say that a single trace $\tau\in S(G,u)$ is \emph{refinable,
weakly good, or good} if the singleton $\{\tau\}$ has the
corresponding property.

There are generalizations, variations, and implications
discussed in \cite{BeH}.  For example, \cite[Lemma 1.1(b)]{BeH}
shows that a trace $\tau\in S(G,u)$ is good if and only
if it is weakly good.  Moreover, a good trace is refinable.
Here we combine our previous results with
\cite[Proposition 7.6]{BeH} to prove the following.

\begin{thm}
\label{ref:thm}
Let $(G,u)$ be a simple dimension group for which $\Aff S(G,u)$ is
finite-dimensional, let $\tau\in S(G,u)$, and let $Z:=Z_G(\ker\tau)
=\Set{\sigma\in S(G,u)}{\ker\tau\subseteq\ker\sigma}$.  The following
conditions are equivalent:
\begin{itemize}
\item[(i)] $\tau$ is refinable,
\item[(ii)] $Z$ is refinable,
\item[(iii)] $Z$ is good,
\item[(iv)] $Z$ is weakly good.
\end{itemize}
When they are met, $G/\ker\tau$ is a simple dimension group;
in particular, it is unperforated.
\end{thm}

\begin{proof}
(i) $\Leftrightarrow$ (ii):
This is because elements $a_1,a_2,b$ of $G$ satisfy
$\tau(a_1)+\tau(a_2)=\tau(b)$
if and only if $\sigma(a_1)+\sigma(a_2)=\sigma(b)$
for each $\sigma\in Z$; similarly $a,a'\in G$ satisfy
$\tau(a)=\tau(a')$ if and only if $\sigma(a)=\sigma(a')$
for each $\sigma\in Z$.

(ii) $\Rightarrow$ (iii):
Suppose that $Z$ is refinable.  Put
$H=\ker\tau=\cap_{\sigma\in Z}\ker\sigma$.
Then, by \cite[Proposition B.5]{BeH}, every trace in $S(G,u)$
maps $H$ to $\{0\}$ or to a dense subgroup of $\bR$.
As $\Aff S(G,u)$ has finite dimension, it follows from
Theorem \ref{suf:thm} that $H$ has property $(B)$ inside $G$ and so,
by Corollary \ref{unperf:cor}, the torsion-free group $G/H$ is
unperforated.  By \cite[Proposition 7.6(a)]{BeH}, this quotient
also has the Riesz interpolation property.  So, $G/H$ is a
simple dimension group. To conclude that $Z$ is a good set
of traces, we need to modify slightly the argument
of \cite[Proposition 7.6(f)]{BeH} as follows.

Let $a\in G$ and $b\in G^+$ with $0<\sigma(a)<\sigma(b)$
for each $\sigma\in Z$.  As the traces
in $S(G/H,u+H)$ are induced by the elements of $Z_G(H)=Z$,
and as $G/H$ is unperforated, we deduce that $a+H$ and
$b-a+H$ belong to $(G/H)^{++}$ \cite[Corollary 4.13]{G}.
Thus, these classes contain elements $a_1$ and $a_2$ of $G^+$,
respectively.  We have $\sigma(a_1)+\sigma(a_2)=\sigma(b)$
for each $\sigma\in Z$.  Since $Z$ is refinable, there
exist $a_1',a_2'\in G^+$ such that $a_1'+a_2'=b$ and
$\sigma(a_i')=\sigma(a_i)$ for each $i=1,2$ and each
$\sigma\in Z$.  In particular, $a_1'\in G^+$ satisfies
$a_1'\le b$ and $\sigma(a_1')=\sigma(a)$ for each
$\sigma\in Z$.  We have $a_1'\neq 0$ and $a_1'\neq b$
since $0<\tau(a_1')<\tau(b)$, thus $0<a_1'<b$ since $G$
is simple.  Therefore $Z$ is good.

(iv) $\Rightarrow$ (ii): Suppose that $Z$ is weakly good, and
that $a_1,a_2,b\in G^+$ satisfy $\sigma(a_1)+\sigma(a_2)=\sigma(b)$
for all $\sigma\in Z$.  If $a_2\neq 0$, then $a_2\in G^{++}$, so
for all $\sigma\in Z$ we have $\sigma(a_1)=\sigma(b)-\sigma(a_2)
< \sigma(b)$. As $Z$ is weakly good, there exists $a_1'\in G^+$
such that $a_1'\le b$ and $\sigma(a_1')=\sigma(a_1)$ for all
$\sigma\in Z$.  Put $a_2'=b-a_1'$.  Then $a_1',a_2'\in G^+$
satisfy $a_1'+a_2'=b$ and $\sigma(a_i')=\sigma(a_i)$ for all $i=1,2$
and all $\sigma\in Z$.  If $a_2=0$, the same holds for the choice
of $a_1'=b$ and $a_2'=0$. Thus, $Z$ is refinable.

This completes the proof since the implication (iii)
$\Rightarrow$ (iv) is immediate and the last
assertion has been established in the course
of proving (ii) $\Rightarrow$ (iii).
\end{proof}

\begin{cor}
Under the same hypotheses, the following conditions are
equivalent:
\begin{itemize}
\item[(i)] $\tau$ is good;
\item[(ii)] $\tau$ is refinable and $Z_G(\ker \tau) = \{\tau\}$.
\end{itemize}
\end{cor}

\begin{proof}
Suppose first that $\tau$ is good, and let
$\sigma\in Z_G(\ker\tau)$.
Since $\ker\tau\subseteq\ker\sigma$, there exists a group
homomorphism $\psi\colon\tau(G)\to \bR$ such that
$\sigma(a)=\psi(\tau(a))$ for each $a\in G$.  We have
$\psi(1)=1$ since $\sigma(u)=\tau(u)=1$.  Moreover,
$\psi$ is order preserving: if $\tau(a)\ge 0$ for
some $a\in G$, then, since $\tau$ is good and $a\le nu$
for some $n\in\bN$, there exists $a'\in G^+$ with
$\tau(a')=\tau(a)$, and so $\psi(\tau(a))=\sigma(a)
=\sigma(a')\ge 0$.  Thus, $\psi$ is the inclusion
of $\tau(G)$ in $\bR$ and therefore $\sigma=\tau$.  As
$\tau$ is refinable, this proves that (i) implies (ii).
The converse follows from Theorem \ref{ref:thm}.
\end{proof}

Suppose now that $(G,u)$ is an arbitrary simple dimension
group $G$ with order unit, and let $\phi\colon G\to\Aff S(G,u)$
be the natural map.  It is shown in \cite[Corollary 1.8]{BeH}
that a trace $\tau\in S(G,u)$ is good if and only if
$\phi(\ker \tau)$ is dense in
$\Set{h \in \Aff S(G,u)}{h(\tau) = 0}$,
and this characterization is very useful.  There is a similar
necessary condition for $\tau$ to be refinable
\cite[Proposition 7.7(e)]{BeH}, but is far from sufficient,
even when $\Aff S(G,u)$ has finite dimension.
However, when $G$ is $\bR^n$ with the strict ordering,
\cite[Appendix 2]{H2} provides a simple geometric description
of the good subsets of $S(\bR^n,\un)=K_n$ of the form
$K_n\cap V$ for a subspace $V$ of $(\bR^n)^*$ (conjectured in
\cite[p.~6295]{BeH}).  Together with Theorem \ref{ref:thm},
this yields a geometric description of the refinable traces of
$\bR^n$.

As shown in \cite[Lemma 7.3]{BeH}, a sufficient condition
for a trace $\tau\in S(G,u)$ to be refinable is that
$\ker \tau=\Inf G$, where $\Inf G=\ker\phi$ is the so-called
\emph{infinitesimal subgroup} of $G$.  By \cite[Proposition 1.7]{GHH},
the collection of such traces is a dense $G_{\delta}$ of $S(G,u)$.
It is contained in the set of refinable traces of $G$.
The next proposition provides a case of equality.

\begin{prop}
Let $G$ be a dense subgroup of $\bR^n$, free of rank $n+1$,
equipped with the strict ordering inherited from $\bR^n$,
and let $\tau\in S(G,u)$ for some $u\in G^{++}$. Then
$\tau$ is refinable if and only if $\ker \tau = \{0\}$.
\end{prop}

A partially ordered abelian group $G$ as in the statement
of the proposition is called a \emph{critical dimension
group of rank $n+1$} (cf., \cite{H1}).

\begin{proof}
Suppose that $\tau$ is refinable, and let $H = \ker \tau$.
Since $H$ is a proper subgroup of $G$, it is discrete in
$\bR^n$.  Let $S$ denote the set of all linear forms
$\sigma\colon\bR^n\to\bR$ for which $\sigma(H)$ is a nonzero
discrete subgroup of $\bR$.  If $H\neq\{0\}$, then $S$ is
a dense subset of $(\bR^n)^*$, stable under multiplication
by positive real numbers; in particular $S$ contains an
element of $S(G,u)$, in contradiction with
\cite[Proposition B.5]{BeH}.  Thus $H$ must be $\{0\}$.
The converse is clear.
\end{proof}

%
%

\section{A class of examples}
\label{sec:ex}

Recall that Theorem B (stated in the introduction) follows
from Theorems \ref{nec:thm} and \ref{suf:thm}.  In this
section, we apply the result to determine, among a class
of subgroups of $\bR^n$, those that have the one-sided
approximation property $(B)$.  We will also use the
following consequence of Theorem B which stresses the
link between properties $(A)$ and $(B)$.

\begin{thm}
\label{ex:thm}
Let $H$ be a subgroup of $\bR^n$, let $F$ be the
smallest face of $K_n$ containing
$Z(H)=\Set{\tau \in K_n}{\tau(H) = \{0\}}$,
let $\tau_{i_1},\dots,\tau_{i_\ell}$ be the vertices
of $F$, and let $\{\tau_{j_1},\dots,\tau_{j_k}\}$
be a maximal set of vertices of $F$ whose restrictions to
$\bR H$ are linearly independent over $\bR$.  Then the
following conditions are equivalent:
\begin{itemize}
\item[(i)] the group $H$ has property $(B)$ in $\bR^n$,
\item[(ii)] its projection
  $H':=\Set{(\tau_{i_1}(h),\dots,\tau_{i_\ell}(h))}{h\in H}$
  has property $(A)$ in $\bR^\ell$,
\item[(iii)] its projection
  $H'':=\Set{(\tau_{j_1}(h),\dots,\tau_{j_k}(h))}{h\in H}$
  is dense in $\bR^k$.
\end{itemize}
\end{thm}

Recall that $\{\tau_1,\dots,\tau_n\}$ denotes the basis of
$(\bR^n)^*$ dual to the canonical basis $\{e_1,\dots,e_n\}$
of $\bR^n$.

\begin{proof}
By Theorem B, condition (i) is equivalent to requiring that
$\tau(H)$ be $\{0\}$ or dense in $\bR$ for each $\tau\in F$,
while by Theorem A, condition (ii) is equivalent to asking
that $\phi(H)$ is $\{0\}$ or dense in $\bR$ for each
$\phi\in \bR F$.  Thus the latter implies the former. To prove
the converse, suppose that condition (i) holds and let
$\phi\in\bR F$.  We have
$\phi=\phi(e_{i_1})\tau_{i_1}+\cdots+\phi(e_{i_\ell})\tau_{i_\ell}$.
By hypothesis, each $\nu\in Z(H)$ admits a similar decomposition
with coefficients $\nu(e_{i_j})\ge 0$ ($1\le j\le \ell$) of
sum $1$ and, for each $j=1,\dots,\ell$, there is at least one
element $\nu_j$ of $Z(H)$ with $\nu_j(e_{i_j})>0$.  Then,
$\nu=\ell^{-1}(\nu_1+\cdots+\nu_\ell)\in Z(H)$ has
$\nu(e_{i_j})>0$ for $j=1,\dots,\ell$. In other words, $\nu$
belongs to the relative interior of $F$.  Choose $a>0$ such
that $c_j:=\phi(e_{i_j})+a\nu(e_{i_j})>0$ for $j=1,\dots,\ell$
and let $c=c_1+\cdots+c_\ell$.  Then $\tau:=c^{-1}(\phi+a\nu)$
belongs to $F$ and so $\tau(H)$ is zero or dense in $\bR$.
Therefore, the same applies to $\phi(H)=c\tau(H)$, showing that
condition (ii) holds.

Finally, the projection map $\pi\colon\bR^\ell\to\bR^k$ given by
$\pi(x_{i_1},\dots,x_{i_\ell})=(x_{j_1},\dots,x_{j_k})$ maps
$H'$ to $H''$ and induces an isomorphism of vector spaces from $\bR H'$
to $\bR H''=\bR^k$.  By Kronecker's theorem, condition (ii) is
equivalent to $H'$ being dense in $\bR H'$.  Thus it is equivalent
to $H''$ being dense in $\bR^k$, which is condition (iii).
\end{proof}

Suppose that $n\ge 2$ and let $\tau=(\tau_1+\tau_2)/2\in K_n$.
We consider subgroups $H$ of $\bR^n$ contained in $\ker\tau$,
or equivalently, for which $\tau\in Z(H)$.  We first note the
following simple consequence of the result above.

\begin{cor}
\label{ex:cor}
With the notation above, suppose that $Z(H)=\{\tau\}$.  Then $H$ has
property $(B)$ if and only if $\tau_1(H)$ is dense in $\bR$.
\end{cor}

\begin{proof}
Here the smallest face $F$ of $K_n$ that contains $Z(H)$ is the
convex hull of $\tau_1$ and $\tau_2$.  As the restriction of $\tau_1$
to $\bR H$ is nonzero while $\tau_2$ coincides with $-\tau_1$ on
$\bR H$, Theorem \ref{ex:thm} applies with $k=1$ and $j_1=1$, and the
conclusion follows.
\end{proof}

\medskip
We now turn to a more specific example.

\begin{example}
\label{ex:example}
Suppose that $n\ge 3$ and let $\alpha_1,\dots,\alpha_{n-1},
\eta_1,\dots,\eta_{n-2}\in\bR$.
Consider the subgroup $H$ of $\bR^n$ generated by the
rows $h_1,\dots,h_{n-1}$ of the matrix
\[
C = \begin{pmatrix} \alpha_1 &-\alpha_1 & &  & &\\
\vdots&\vdots  &&& I_{n-2} && \\\
\alpha_{n-2} &-\alpha_{n-2}&  &    &\\
\alpha_{n-1} &-\alpha_{n-1} & \eta_1 & \eta_2 &\dots & \eta_{n-2}
\end{pmatrix},
\]
where $I_{n-2}$ is the identity matrix of size $n-2$.
Let $S$ denote the set of indices $j$ with $1\le j\le n-2$
for which $\alpha_j\neq 0$.
Then $H$ has property $(B)$ if and only if
$\tau_1(H)=\bZ\alpha_1+\cdots+\bZ\alpha_{n-1}$ is zero or
dense in $\bR$ and one of the following mutually exclusive
conditions holds.
\begin{itemize}
\item[(i)] The rank of $C$ is $n-1$.
\item[(ii)] The rank of $C$ is $n-2$ and not all $\alpha_j$
with $j\in S$ have the same sign.
\item[(iii)] The rank of $C$ is $n-2$, all $\alpha_j$
with $j\in S$ have the same sign, and the set
$\{1\}\cup\Set{\eta_j}{j\in S}$ is linearly independent
over $\bQ$.
\end{itemize}
\end{example}

\medskip
As the proof will show, we have $Z(H)=\{\tau\}$ in cases
(i) and (ii), while $Z(H)$ is a line segment in case (iii).

\begin{proof}
Since $\tau\in Z(H)$ and since $\tau_1$ belongs to the smallest
face of $K_n$ containing $\tau$, the condition that $\tau_1(H)$
is zero or dense in $\bR$ is necessary by Theorem B.
Suppose that it holds.
If the rank of $C$ is $n-1$, then $\bR\,H = \ker\tau$,
thus $Z(H)=K_n\cap\bR\tau=\{\tau\}$
and so $H$ has property $(B)$ by Corollary \ref{ex:cor}.
Otherwise, the rank of $C$ is $n-2$ and we have
\[
 h_{n-1}=\eta_1h_1+\cdots+\eta_{n-2}h_{n-2}.
\]
Moreover, the linear form $\phi$ given by
\[
 \phi = \tau_2 + \alpha_1\tau_3+\cdots+\alpha_{n-2}\tau_n
      = \tau_2+\sum_{j \in S} \alpha_j \tau_{j+2}
\]
vanishes on $h_1,\dots,h_{n-2}$ and therefore on the whole
of $H$.  This means that the annihilator of $H$ in $(\bR^n)^*$
is $\bR\tau+\bR\phi$ and so
\[
 Z(H) = K_n \cap (\bR\tau+\bR\phi).
\]
If the $\alpha_j$ with $j\in S$ do not all have the same sign,
this implies that $Z(H)=K_n \cap\bR\tau = \{\tau\}$ and again
$H$ has property $(B)$ by Corollary \ref{ex:cor}.
Otherwise $\phi$ or $2\tau-\phi$
is a positive linear functional and we find that $Z(H)$ is
the line segment in $K_n$ joining $\tau$ to either $c^{-1}\phi$ or
$c^{-1}(2\tau-\phi)$, where $c=1+\sum_{j\in S}|\alpha_j|$.
Then the smallest face $F$ of $K_n$ containing $Z(H)$ is the
convex hull of the set
$\{\tau_1,\tau_2\}\cup\Set{\tau_{j+2}}{j\in S}$ and
$\Set{\tau_{j+2}}{j\in S}$ is a maximal set of vertices of $F$
whose restriction to $\bR H$ are linearly independent over $\bR$.
Letting $j_1,\dots,j_k$ denote the distinct elements of $S$,
we deduce from Theorem \ref{ex:thm}, that $H$ has property $(B)$
if and only if its projection
\[
 \Set{(t_{j_1+2}(h),\dots,\tau_{j_k+2}(h))}{h\in H}
 = \bZ^k+\bZ(\eta_{j_1},\dots,\eta_{j_k}).
\]
is dense in $\bR^k$, that is if and only if $1,\eta_{j_1},\dots,
\eta_{j_k}$ are linearly independent over $\bQ$.
\end{proof}

\appendix

%
%
%

\section{Gordan's theorem and Farkas' lemma}
\label{sec:A}

Gordan's theorem (not to be confused with Gordan's
lemma, a result--by the same Gordan-- concerning toric
varieties) asserts the following \cite{Gor}.

\begin{thm}[Gordan] Let $A$ be an $m \times n$ real matrix. Exactly
one of the following is true.
\begin{itemize}
\item[(i)] There exists $\uy \in (\bR^n)^+\setminus \{\pmb 0\}$
such that $A\uy^t = \uzero$.
\item[(ii)] There exists $\ux \in \bR^m$ such that $\ux A\in (\bR^m)^{++}$.
\end{itemize}
\end{thm}

Here, we view the elements of $\bR^m$ and of $\bR^n$ as row vectors.
If (ii) holds, then we may choose $\ux$ in $\bZ^m$, by a
simple density argument.  Thus Gordan's theorem is the special case
of Proposition \ref{suf:prop:ZHempty} applied to the group $G = \bR^n$
with the usual coordinatewise ordering, and to the subgroup $H$
of $G$ generated by the rows of $A$:  alternative (i) says that
there exists a trace $\tau$ of $G$ such that $\tau(H)=\{0\}$,
while alternative (ii) with $\ux \in \bZ^m$ says that $H$
contains an element of $G^{++}$.

Gordan's theorem was followed chronologically by Farkas' lemma \cite{F},
which is now typically  used to prove the former.  We state Farkas' Lemma
below in one of its numerous equivalent forms.  It can be used to
provide a direct proof of Theorem B.

\begin{thm}[Farkas]
\label{app:Farkas}
Let $A$ be a real $m \times n$ matrix, and let
$\ub \in \bR^n$. Exactly one of the following is true.
\begin{itemize}
\item[(i)] There exists  $\uy \in (\bR^n)^+$ such
  that $A\uy^t = \uzero$ and $\ub\uy^t< 0$.
\item[(ii)] There exists $\ux \in \bR^m$ such that
  $\ux A \leq \ub$.
\end{itemize}
\end{thm}

Contrary to Gordan's theorem, this result has limited
extension to partially ordered abelian groups.
To explain, let $G = \bR^n$ equipped
with the coordinatewise ordering, and let $H$
be the subspace of $\bR^n$ generated, over $\bR$, by the
rows of $A$, where $A$ is as in the statement of
Theorem \ref{app:Farkas}.  Then the result says that,
for any $g\in G$, either (ii) there exists $h\in H$
such that $h\le g$, or (i) there exists a trace
$\tau\in S(G,\un)$ such that $\tau(H)=\{0\}$ and
$\tau(g)<0$. Equivalently, this means that, for any
$g\in G$, we have
\begin{equation}
\label{app:cond}
 (g+H)\cap G^+\neq \{0\}
 \ \Longleftrightarrow\
 \tau(g)\ge 0 \text{ for all } \tau\in Z_G(H).
\end{equation}
This property is easy to characterize when $H\cap G^+=\{0\}$,
or more generally, when $H$ is a convex subgroup of $G$.

\begin{lemma}
\label{app:lemma}
Let $(G,u)$ be a partially ordered abelian group with
order unit, and let $H$ be a convex subgroup of $G$.
Then \eqref{app:cond} holds for each $g\in G$ if and only if
the partially ordered group $G/H$ is archimedean.
\end{lemma}

Recall that a partially ordered abelian group $G$ is
{\it archimedean\/} if, for $x,y \in G$, the condition $nx \leq y$
for all positive integers $n$ entails that $-x \in G^+$
(a inequivalent definition, frequently seen, requires
$x \in G^+$ at the outset). When $G$ has an order unit $u$,
it is archimedean if and only if the traces in $S(G,u)$
determine the ordering on $G$, that is,
$G^+=\Set{g\in G}{\hat g \geq 0}$
(\cite[Theorem 4.14]{G}).  Under the hypotheses of
Lemma \ref{app:lemma}, this is exactly
what \eqref{app:cond} means for the pair $(G/H,u+H)$.

Archimedeanness is a strong property, particularly if the partially
ordered abelian group is simple.  It was not known until 2013 that for
every (infinite-dimensional) metrizable Choquet simplex $K$, there
exists a simple archimedean dimension group $(G,u)$ whose trace
space is $K$ \cite{H3}.

\enddocument